\input amstex
\documentstyle{amsppt}
\NoPageNumbers		

\mag=1200
\addto\tenpoint{\normalbaselineskip=15pt\normalbaselines}	
\hsize=140 true mm      \hoffset=8 true mm
\vsize=216 true mm      \voffset=7 true mm
\relpenalty=10000
\binoppenalty=10000
\parskip=3pt plus1pt
\font\titlefont=cmss12 scaled 1095

\catcode`\@=11
\newfam\ssffam
\font@\tenssf=cmss10
\font@\eightssf=cmss8
\addto\tenpoint{%
  \textfont\ssffam=\tenssf \scriptfont\ssffam=\eightssf}%
\addto\eightpoint{%
  \textfont\ssffam=\eightssf \scriptfont\ssffam=\eightssf}%
\def\proclaimheadfont@{\smc} \def\headfont@{\tenssf}
\catcode`\@=\active

\hrule height0pt
\vskip15mm
\topmatter

\title\nofrills\titlefont
Writing representations over minimal fields
\endtitle


\leftheadtext\nofrills{S.\,P. Glasby and R.\,B. Howlett}
\rightheadtext\nofrills{Writing representations over minimal fields}

\author
\vskip-4mm
\tenssf S.\,P. Glasby and R.\,B. Howlett
\vskip2mm
{\eightssf School of Mathematics and Statistics\\
University of Sydney, NSW 2006, Australia}\\
\vskip3mm
\endauthor

\def\E{\Bbb E}
\def\F{\Bbb F}
\def\K{\Bbb K}
\def\GL{\text{GL}}
\def\Mat{\text{Mat}}
\def\refHR{1}
\def\refHB{2}
\def\qed{\hfill\llap{$\sqcup$}\llap{$\sqcap$}}

\abstract\baselineskip12pt
The chief aim of this paper is to describe a procedure
which, given a $d$-dimensional absolutely irreducible matrix representation 
of a finite group over a finite field $\E$, produces an equivalent
representation such that all matrix entries lie in a subfield $\F$ 
of $\E$ which is as small as possible. The algorithm relies on a matrix
version of Hilbert's Theorem~90, and is probabilistic with expected running
time $O(|\E:\F|d^3)$ when $|\F|$ is bounded. Using similar methods
we then describe an algorithm which takes as input a prime number and a
power-conjugate presentation for a finite soluble group, and as output
produces a full set of absolutely irreducible representations of the 
group over fields whose characteristic is the specified prime, each
representation being written over its minimal field. 
\endabstract


\endtopmatter

\document

\hrule height0pt
\vskip-3mm
\head
1. The main algorithm
\endhead

\noindent Let $\rho\colon G\to\GL(d,\E)$ be an absolutely irreducible
representation of the group $G$. It is clear that there exists a
subfield $\F$ of $\E$, minimal
with respect to inclusion, such that there exists a representation
$G\to\GL(d,\F)$ equivalent to $\rho$. If $\E$ has nonzero
characteristic, then $\F$ is determined by $\rho$, and coincides with
the subfield generated by the character values of $\rho$
(see [\refHB, VII Theorem 1.17]). Indeed, the arguments presented here yield
a proof of this fact. If $\E$ has characteristic zero,
there may be more than one choice for $\F$.

Suppose that $\F$ is a subfield of $\E$ such that $\E$ is a finite
Galois extension of $\F$ whose Galois group is cyclic,
of order~$t$, and generated by~$\alpha$. Assume further that the norm map
from $\E$ to  $\F$ (given by $\lambda\mapsto
\lambda\lambda^\alpha \lambda^{\alpha^2}\cdots\lambda^{\alpha^{t-1}}$)
is surjective. This hypothesis certainly holds if $|\E|$ is finite,
and this is the case of principal interest to us.
Our first objective is to describe a procedure which determines whether
an absolutely irreducible representation $\rho\colon G\to
\GL(d,\E)$ of a finitely generated group $G$ is equivalent
to a representation $G\to\GL(d,\F)$, and if so, finds an
$A\in\GL(d,\E)$ such that $A^{-1}\rho(g)A\in\GL(d,\F)$
for all $g\in G$. Note that if $g_1,\,g_2,\,\ldots,\,g_n$ generate $G$,
this condition is equivalent to $A^{-1}\rho(g_i)A\in\GL(d,\F)$
for all $i\in\{1,2,\ldots,n\}$.

A basic step in our algorithm involves testing whether two given matrix
representations of $G$ are equivalent, and if they
are, finding a nonsingular intertwining matrix. The naive approach to 
this problem involves solving $nd^2$ homogeneous linear equations in 
$d^2$ unknowns over the field $\E$. Computationally, this has cost 
$O(nd^6)$. Alternatively, there is a probabilistic algorithm, described 
by Holt and Rees in [\refHR], which has expected running time $O(d^3)$.
(This complexity result, and those throughout this section, assume that the cost
of field arithmetic, including applying a field automorphism, is $O(1)$.)

With the notation as above, suppose that $A\in\GL(d,\E)$ 
has the property that $A^{-1}\rho(g)A\in\GL(d,\F)$ for all~$g\in G$. The
automorphism $\alpha$ of $\E$ gives rise to an automorphism of
$\Mat(d,\E)$ (the algebra of $d\times d$ matrices
over~$\E$) which we also denote by~$\alpha$. Since the fixed subfield
of $\alpha$ is $\F$, it is clear that $B\in\Mat(d,\E)$
satisfies $B^\alpha=B$ if and only if $B\in\Mat(d,\F)$. So
$(A^{-1}\rho(g)A)^\alpha=A^{-1}\rho(g)A$ for all~$g\in G$, and thus
$C=A(A^\alpha)^{-1}$ satisfies
$$
C^{-1}\rho(g)C=\rho(g)^\alpha \qquad\text{(for all~$g\in G$).} \leqno(1)
$$
Since $\rho$ is absolutely irreducible, equation (1) determines $C$ 
up to a nonzero scalar multiple. The first step in our procedure is,
therefore,
to use an algorithm such as in [\refHR] to find (if possible) a
$C\in\GL(d,\E)$ satisfying~(1). If no such $C$ exists, then
$\rho$ cannot be written over $\F$; so assume henceforth that such a
$C$ has been found.

\proclaim{Proposition (1.1)} If $C\in\GL(d,\E)$ satisfies~$(1)$, then
$CC^\alpha C^{\alpha^2}\cdots C^{\alpha^{t-1}}$ equals $\mu I$ where
$\mu\in \F$ and $I$ is the $d\times d$ identity matrix.
\endproclaim

\demo{Proof} Since
$CC^\alpha C^{\alpha^2}\cdots C^{\alpha^{t-1}}\!$ conjugates $\rho(g)$ 
to $\rho(g)^{\alpha^t}\!=\rho(g)$ for all~$g\in G$,
it must equal $\mu I$ for some $\mu\in\E$, since $\rho$ is assumed 
to be absolutely irreducible. However,
$$
\mu^\alpha I=C(\mu I)^\alpha C^{-1}=
C(C^\alpha C^{\alpha^2}C^{\alpha^3}\cdots C^{\alpha^t})C^{-1}
=CC^\alpha C^{\alpha^2}\cdots C^{\alpha^{t-1}}=\mu I,
$$
and so $\mu\in\F$, as desired.\qed
\enddemo

The computation of $\mu$ can be effected by $t-1$ vector by matrix
multiplications, since if $v$ is the first row of~$C$ then $\mu$ is the
first component of the row vector
$vC^\alpha C^{\alpha^2}\cdots C^{\alpha^{t-1}}$. This has cost 
$O(td^2)$. If $t$ is large compared with~$d$, then $\mu$ may be computed
at cost $O((\log t)d^3)$ by using the fact that
$C_{2i}=C_i(C_i)^{\alpha^i}$ for each $i$, where
$C_i=CC^\alpha \cdots C^{\alpha^{i-1}}$.

Since the norm map from $\E$ to $\F$ is assumed to be 
surjective, there exists a $\nu\in\E$ whose norm is $\mu$. We 
do not address here the practical problem of finding $\nu$ 
given~$\mu$. The methods used for storing field elements and performing
field computations obviously affect this issue. (When $|\F|$ is bounded,
there is an $O(1)$ probalistic algorithm for computing $\nu$.)
Once $\nu$ has been
found we may replace $C$ by $\nu^{-1}C$, and assume thereafter that
$CC^\alpha\cdots C^{\alpha^{t-1}}=I$. 

\proclaim{Lemma (1.2)} If $C\in\GL(d,\E)$ satisfies
$CC^\alpha\cdots C^{\alpha^{t-1}}=I$, then there exists a nonzero
column vector $v\in\E^d$ such that $Cv^\alpha=v$.
\endproclaim

\demo{Proof} Let $u_0\in\E^d$ be nonzero, and for $i>0$ define
$u_i$ recursively by $u_i=Cu_{i-1}^\alpha$. Observe that $u_t=u_0$. Now 
since the field automorphisms $\alpha^0$, $\alpha^1$, \dots,
$\alpha^{t-1}$ are distinct they are linearly independent, and since the
$u_i$ are nonzero it follows that there exists a $\lambda\in\E$ such that 
$v=\sum_{i=0}^{t-1}\lambda^{\alpha^i}u_i\ne 0$. Moreover, 
$Cv^\alpha=\sum_{i=1}^t\lambda^{\alpha^i}Cu_{i-1}^\alpha=v$, as desired.\qed
\enddemo

The following proposition may be viewed as a generalization of the 
multiplicative form of Hilbert's Theorem~90. The corresponding 
generalization of the additive form is trivially true.

\proclaim{Proposition (1.3)} If $C\in\GL(d,\E)$ satisfies
$CC^\alpha\cdots C^{\alpha^{t-1}}=I$, then there exists an 
$A\in\GL(d,\E)$ with $C=A(A^\alpha)^{-1}$.
\endproclaim

\demo{Proof} The result is true when $d=1$ by the multiplicative form of
Hilbert's Theorem~90. Proceeding by induction, assume that $d>1$.
By Lemma~(1.2) there exists a nonzero vector~$v$ such that 
$Cv^\alpha=v$, and if $B$ is an invertible matrix with $v$ as its first 
column then
$$
B^{-1}CB^\alpha=\pmatrix 1&u\\0&C_1\endpmatrix
$$
where $C_1\in\GL(d-1,\E)$ satisfies
$C_1C_1^\alpha\cdots C_1^{\alpha^{t-1}}=I$. By the inductive hypothesis,
there exists an $A_1\in\GL(d-1,\E)$ such that
$C_1=A_1(A_1^\alpha)^{-1}$, and it follows that
$$
\pmatrix 1&0\\0&A_1\endpmatrix^{-1}B^{-1}CB^\alpha
\pmatrix 1&0\\0&A_1\endpmatrix^\alpha
=\pmatrix 1&u_1\\0&I\endpmatrix
$$
where $u_1=u(A_1^{-1})^\alpha$ satisfies $\sum_{i=0}^{t-1}u_1^{\alpha^i}=0$.
It follows from the additive form of Hilbert's Theorem~90 that
there exists a row vector $u_2$ with $u_1=u_2-u_2^\alpha$, and then
$$
A=B\pmatrix 1&0\\0&A_1\endpmatrix\pmatrix 1&u_2\\0&I\endpmatrix
$$
has the required property $C=A(A^\alpha)^{-1}$.\qed
\enddemo

Note that if $C=A(A^\alpha)^{-1}$ then the map $\Mat(d,\E)\to\Mat(d,\E)$
given by
$$
\align
X&\mapsto X+CX^\alpha+CC^\alpha X^{\alpha^2}+\cdots+CC^\alpha\cdots 
C^{\alpha^{t-2}}X^{\alpha^{t-1}}\\
&=A(A^{-1}X+(A^{-1}X)^\alpha+\cdots+(A^{-1}X)^{\alpha^{t-1}})
\endalign
$$
has image consisting of all matrices of the form $AY$ with
$Y\in\Mat(d,\F)$. These are exactly the matrices $A'\in\Mat(d,\E)$ 
such that $(A^{-1}A')^\alpha=A^{-1}A'$, or equivalently,
$C(A')^\alpha=A'$. If $X$ is chosen arbitrarily and $X\mapsto AY=A'$,
then the probability 
that $Y$ is invertible (so that $C=A'((A')^\alpha)^{-1}$) is
$|\GL(d,\F)|/|\Mat(d,\F)|$. It follows that a reasonable procedure for
finding an $A$ satisfying the equation $C=A(A^\alpha)^{-1}$ is to choose
$X\in\Mat(d,\E)$ randomly and compute
$A=X+CX^\alpha+CC^\alpha X^{\alpha^2}+\cdots+CC^\alpha\cdots 
C^{\alpha^{t-2}}X^{\alpha^{t-1}}$, repeating if necessary until an 
invertible~$A$ is found. (One may show that
$1-|\F|^{-1}\ge|\GL(d,\F)|/|\Mat(d,\F)| > 1-|\F|^{-1}-|\F|^{-2}\ge 1/4$.)

Observe that $C=A(A^\alpha)^{-1}$ combines with equation~(1) to give
$$
A^{-1}\rho(g)A=(A^{-1}\rho(g)A)^\alpha \qquad\text{(for all~$g\in G$).}
$$
It follows that $A^{-1}\rho(g)A\in\GL(d,\F)$ for each~$g$,
and we have achieved our goal of constructing a representation 
equivalent to $\rho$ with image contained in $\GL(d,\F)$.
Note that if $A_i=X+CX^\alpha+CC^\alpha X^{\alpha^2}
+\cdots+CC^\alpha\cdots C^{\alpha^{i-2}}X^{\alpha^{i-1}}$ then
$A_{i+1}=X+CA_i^\alpha$, and it follows that $A_t$ can be evaluated with
$t-1$ matrix multiplications and $t-1$ matrix additions. It can be seen,
therefore, that our procedure has expected running time
$O(|\E:\F|d^3)$.

\head
2. Absolutely irreducible representations of soluble groups
\endhead

\noindent Suppose that we are given a consistent power-conjugate presentation for
a finite group~$G$. That is, $G$ is generated by 
$g_1,\,g_2,\,\ldots,\,g_n$, where $n$ is the composition length of $G$,
with defining relations
$$
\alignat2
g_i^{p_i}&=v_i\qquad&&\text{($1\le i\le n$)}\\
g_i^{-1}g_jg_i&=w_{ij}&&\text{($1\le i<j\le n$)}
\endalignat
$$
where each $p_i$ is a prime and each $v_i$ is a word in the generators
$g_j$ for $i<j\le n$, and each $w_{ij}$ is a word in the $g_k$ for
$i<k\le n$. It is clear that a group has such a presentation if and only
if it is finite and soluble. Specifically, if $G_i$ is the subgroup of $G$ 
generated by $g_i,g_{i+1},\ldots,g_n$, then
$$
G=G_1\ge G_2\ge \cdots \ge G_n\ge G_{n+1}=\{1\}\leqno({*})
$$
is a subnormal series, and for each~$i$ the quotient $G_i/G_{i+1}$ has 
order dividing~$p_i$. Given that $n$ is the composition length of $G$,
it follows that $({*})$ is a composition 
series and the order of $G_i/G_{i+1}$ is exactly~$p_i$. We will show how
the natural algorithm for constructing the absolutely irreducible 
representations of~$G$ (in a fixed nonzero characteristic),
by working up the 
composition series~$({*})$, can be readily adapted to ensure that each
representation is written over its minimal field. We consider that we 
have constructed a representation of the group $G_i$ once we have 
computed matrices representing the generators
$g_i$,~$g_{i+1}$,~\dots,~$g_n$.

For ease of exposition we let $\K$ be a fixed algebraic closure of
a field of prime order, and deal henceforth only with subfields of
$\K$. Assume, inductively, that we have constructed representations
$\sigma_1,\,\sigma_2,\,\ldots,\,\sigma_s$ of the group $G_2$
such that
\roster
\item"(i)"
each $\sigma_i$ is absolutely irreducible and written
over its (unique) minimal subfield of $\K$, and
\item"(ii)"
every absolutely irreducible representation of $G_2$ over
$\K$ is equivalent to exactly one of the~$\sigma_i$.
\endroster
Henceforth, to simplify the notation, we write $H=G_2$, $a=g_1$ and $p=p_1$.

The absolutely irreducible $\K$-representations of $H$ are permuted by 
$G$ via
$$\sigma^g(h)=\sigma(ghg^{-1})
$$
for all $h\in H$ and $g\in G$.
The first step is to find, for each $i$, which of the representations
$\sigma_1,\,\sigma_2,\,\ldots,\sigma_s$ is equivalent to
the representation $\sigma_i^a$.
If $\sigma_i^a$ is equivalent to $\sigma_i$, then there exists a 
representation of $G$ extending $\sigma_i$; the minimal field for any such 
extension will be an extension of the field of $\sigma_i$. If 
$\sigma_i^a$ is not equivalent to $\sigma_i$, then $\sigma_i$ will be 
$G$-conjugate to $p=|G:H|$ of the representations $\sigma_k$. In this 
case the representation of $G$ induced from~$\sigma_i$ is absolutely 
irreducible; however, its minimal field may be smaller than that
of~$\sigma_i$.
Since $G$-conjugate representations of $H$ 
yield equivalent induced representations of~$G$, one representative only
should be chosen from each $G$-conjugacy class.

\noindent{\smc Case 1.} {\it Assume that $\E$ is a finite field, and
$\sigma\colon H\to\GL(d,\E)$ is an absolutely irreducible representation,
with minimal field $\E$, such that $\sigma^a$ is equivalent to~$\sigma$.}

Compute a matrix $A\in\GL(d,\E)$ such that
$A\sigma(h)A^{-1}=\sigma(aha^{-1})$ for all $h\in H$.
As $\sigma$ is absolutely irreducible and $a^p\in H$, so
$A^p=\mu\sigma(a^p)$ for some $\mu$ in $\E^\times$ (the multiplicative 
group of $\E$). If the characteristic of $\E$ equals~$p$, then $\mu$ has 
a unique $p$th root $\nu\in\E^\times$. Indeed, $\nu$ is a power of $\mu$
since~$p$ is coprime to $|\E^\times|$.
In this case there is a unique representation $\rho$ of $G$ extending
$\sigma$, given by $\rho(a)=\nu^{-1}A$ and $\rho(h)=\sigma(h)$ for
all $h\in H$. Suppose alternatively that the characteristic
of $\E$ is not $p$. In this case
$\nu^p=\mu$ has exactly $p$ solutions $\nu_1,\ldots,\nu_p$ in $\K$, and
correspondingly there are $p$ pairwise inequivalent extensions 
$\rho_1,\ldots,\rho_p$ of $\sigma$
given by defining $\rho_i(a)=\nu_i^{-1}A$. For each $i$, the
extension field $\E(\nu_i)$ is the minimal field for $\rho_i$.
If $|\E^\times|$ is coprime to~$p$, then one of the
solutions of $\nu^p=\mu$ lies in the field $\E$, while the remaining~$p-1$
solutions generate the same field, which is the 
smallest extension of $\E$ whose order is congruent to~1 modulo~$p$.
If $|\E^\times|$ is a multiple of $p$, then all solutions of $\nu^p=\mu$
generate the same extension $\E'$ of $\E$. Note that $|\E':\E|$ is
$1$ or $p$, and $\E'$ is the smallest extension of $\E$
whose order is congruent to 1~modulo~$p|\nu|$.

\noindent{\smc Case 2.} {\it Assume that $\E$ is a finite field, and 
$\sigma\colon H\to\GL(d,\E)$ is an absolutely irreducible representation,
with minimal field $\E$, such that $\sigma^a$ is not equivalent
to~$\sigma$.}

Let $k$ be the degree of $\E$ over its prime subfield.
If $k$ is not a multiple of~$p$, then $\E$ is the minimal
field for the induced representation $\sigma^G$. 
If $k$ is a multiple of $p$, then $\E$ has an automorphism
$\alpha$ of order $p$ whose fixed subfield,~$\F$, is uniquely defined by
$|\E:\F|=p$.
In this case, if the representation $\sigma^\alpha\colon h\mapsto
\sigma(h)^\alpha$ is not equivalent to one of the $G$-conjugates of
$\sigma$, then $\E$ is still the minimal field for~$\sigma^G$; however,
if $\sigma^\alpha$ is equivalent to 
a $G$-conjugate of $\sigma$ then one can readily show that
$\sigma^G$ is equivalent to $(\sigma^G)^\alpha$, and so the minimal 
field of $\sigma^G$ is $\F$.

We present an explicit construction for an $\F$-representation 
equivalent to $\sigma^G$ in the case that $\sigma^\alpha$ is equivalent 
to a $G$-conjugate of~$\sigma$. Replacing $\alpha$ by a power of itself,
we may assume that $\sigma^\alpha$ is equivalent to~$\sigma^a$. Find an
$A\in\GL(d,\E)$ such that
$$
A\sigma(h)^\alpha A^{-1}=\sigma(aha^{-1}) \qquad
\text{(for all $h\in H$),}\leqno{(2)}
$$
and note that, by absolute irreducibility,
$AA^\alpha\cdots A^{\alpha^{p-1}}=\mu\sigma(a^p)$
for some $\mu\in \E$. As in Proposition (1.1) we see that $\mu\in\F$, 
since
$$
\align
\mu^\alpha\sigma(a^p)^\alpha&=A^\alpha A^{\alpha^2}\cdots A^{\alpha^p}\\
&=A^{-1}(AA^\alpha A^{\alpha^2}\cdots A^{\alpha^{p-1}})A\\
&=\mu(A^{-1}\sigma(a^p)A)\\
&=\mu(A^{-1}\sigma(aa^pa^{-1})A)\\
&=\mu\sigma(a^p)^\alpha,
\endalign
$$
where the last step follows from (2). Hence replacing $A$ by 
$\nu^{-1}A$, where $\nu\in\E^\times$ satisfies
$\nu\nu^\alpha\cdots\nu^{\alpha^{p-1}}=\mu$, we may assume that
$$
AA^\alpha\cdots A^{\alpha^{p-1}}=\sigma(a^p).\leqno{(3)}
$$

The regular representation of $\E$ considered as an $\F$-algebra yields 
an $\F$-algebra monomorphism $\phi\colon\E\to\Mat(p,\F)$, and since 
$\alpha$ is an $\F$-automorphism of $\E$ there is an $M\in\GL(p,\F)$ 
satisfying $M^p=I$ and
$$
M^{-1}\phi(\lambda)M=\phi(\lambda^\alpha)\qquad\text{(for all
$\lambda\in\E$).}
$$
(We remark that computing $\phi$ and $M$ is best done when the elements
of $\E$ are represented as polynomials over $\F$ modulo an irreducible
polynomial. In this case, the assumption in Section~1, that field
arithmetic in $\E$ can be performed in constant time, does not hold.)
Let $\Phi\colon\Mat(d,\E)\to\Mat(pd,\F)$ be defined by
$\Phi((\lambda_{i,j}))=(\phi(\lambda_{i,j}))$, and define $S\in\GL(d,\F)$
to be the diagonal sum of $d$ copies of~$M$. Then $\Phi$ is an
$\F$-algebra monomorphism, and
$$
S^{-1}\Phi(X)S=\Phi(X^\alpha)\qquad\text{(for all $X\in\Mat(d,\E)$).}
\leqno{(4)}
$$
It now follows that there is a representation $\rho\colon G\to\GL(pd,\F)$
such that $\rho(a)=\Phi(A)S^{-1}$ and $\rho(h)=\Phi(\sigma(h))$ for all
$h\in H$, since
$$
\alignat2
\rho(a)^p
&=(\Phi(A)S^{-1})^p\\
&=\Phi(A)(S^{-1}\Phi(A)S)\cdots(S^{-(p-1)}\Phi(A)&&S^{p-1})S^{-p}\\
&=\Phi(A)\Phi(A^\alpha)\cdots\Phi(A^{\alpha^{p-1}})&&
    \text{(using (4) and $S^p=I$)}\\
&=\Phi(\sigma(a^p))&&\text{(by (3))}\\
&=\rho(a^p)
\endalignat
$$
and
$$
\alignat2
\rho(a)\rho(h)\rho(a)^{-1}
&=\Phi(A)S^{-1}\Phi(\sigma(h))S\Phi(A)^{-1}\\
&=\Phi(A)\Phi(\sigma(h)^\alpha)\Phi(A^{-1})\hskip20mm&&\text{(by (4))}\\
&=\Phi(\sigma(aha^{-1}))&&\text{(by (2))}\\
&=\rho(aha^{-1}).\\
\endalignat
$$

It remains to check that $\rho$ is equivalent to~$\sigma^G$.
It is clear that there exists a $T\in\GL(p,\E)$ such that
$$
T\phi(\lambda)T^{-1}=
\text{diag}(\lambda,\lambda^\alpha,\ldots,\lambda^{\alpha^{p-1}})
$$
for all $\lambda\in\E$. Furthermore, if $v_i$ denotes the
$(i+1)$th row of $T$ and
$V_i$ denotes the subspace of $\E^{pd}$ comprising the
elements of the form $(\lambda_1v_i,\lambda_2v_i,\ldots,\lambda_dv_i)$
where $\lambda_1,\,\lambda_2,\,\ldots,\,\lambda_d\in\E$,
then
\roster
\item"(i)" $\E^{pd}= V_0\oplus V_1\oplus\cdots\oplus V_{p-1}$,
\item"(ii)" each $V_i$ is $\rho(H)$-invariant, inducing an action equivalent to 
$\sigma^{a^i}$, and
\item"(iii)" $V_i\rho(a)=V_{i+1}$, where the subscripts are read modulo $p$.
\endroster
Note that (ii) follows from $v_i\phi(\lambda)=\lambda^{\alpha^i}v_i$,
and (iii) follows from the equation $\rho(a)\rho(h)\rho(a)^{-1}=\rho(aha^{-1})$.
These conditions guarantee that $\rho$ is equivalent to
$\sigma^G$, as required.
We have thus achieved our goal of constructing the absolutely 
irreducible representations of $G$ over their minimal fields.

\enddocument
\end